\documentclass[10pt]{article}
\usepackage{amsmath,amssymb}
\usepackage{cite}

\numberwithin{equation}{section}

\newtheorem{remark}{Remark}

\def\pr{\preceq}
\def\be{\begin{equation}}
\def\ee{\end{equation}}
\newcommand{\eps}{\varepsilon}

\author{ Hans-G\"org Roos, TU Dresden}
\title{Layer-adapted meshes for weak boundary layers}
\begin{document}
\maketitle

\section{Second order problems and linear finite elements}

We start with problems of {\it convection-diffusion} type and consider
\begin{equation}
-\eps u''-bu'+cu=f, \quad u(1)=0,\,u'(0)=0,
\end{equation}
assuming V-ellipticity of the related bilinear form, $b>1$ and $0<\eps<<1$.

The boundary condition at $x=0$ implies the existence of a {\it weak} boundary layer, i.e., we have
in a solution decomposition into a smooth part $S$ and a layer part $E$
\[
    |E^{(k)}|\preceq \eps^{1-k}e^{-x/\eps}.
\]
It follows for the Sobolev seminorms
\be |E|_1\pr \eps^{1/2},\quad |E|_2\pr \eps^{-1/2}.
\ee

Define the $\eps$-weighted $H^1$ norm by
\[
   \|v\|_\eps:=\eps^{1/2}|v|_1+|v|_0.
\]
It was already observed in \cite{RR11}, that now the discretization with linear elements on an equidistant mesh
with the mesh size $H$ leads to a uniform (with respect to $\eps$) error estimate in the $\eps$-weighted $H^1$ norm
(for the upwind finite difference scheme, see \cite{Li05}) :
\be
   \|u-u_H\|_\eps \pr H.
\ee
This estimate is easy to prove. First we get for the interpolation error
\[
   \eps^{1/2}|E-E^I|_1\pr \eps^{1/2}H|E|_2\pr H
\]
and
\be\label{1}
  |E-E^I|_0\pr H|E|_1\pr \eps^{1/2}H.
\ee
The factor $\eps^{1/2}$ in \eqref{1} allows to estimate the convective term in the error equation.

Next we study  problems of {\it reaction-diffusion} type
\begin{equation}
-\eps^2 u''+cu=f, \quad u'(0)=0,\,u'(1)=0
\end{equation}
with $c>1$. Again we have weak layers, now at $x=0$ and at $x=1$. In the energy norm
\[
   \|v\|_\eps:=\eps|v|_1+|v|_0
\]
one gets immediately on an equidistant mesh
\be
   \|u-u_H\|_\eps \pr \eps^{1/2}H+H^2.
\ee
One can also obtain an estimate in the balanced norm (see \cite{Ro17} for a survey)
\[
   \|v\|_b:=\eps^{1/2}|v|_1+|v|_0.
\]
Here one uses the $L_2$ projection for the smooth part. On a uniform mesh the projection is
stable in the sense
\[
   |\pi S|_1\pr |S|_1
\]
(see \cite{CT87}).
It follows
\[
  |S-\pi S|_1\pr |S-S^I|_1,
\]
and finally
\be
   \|u-u_H\|_b \pr H.
\ee

\section{Second order problems and higher order finite elements}
For higher order elements the use of a uniform mesh does not lead to uniform convergence. Of course,
it is possible to use the same meshes as for strong layers.

But we try to use coarser meshes, especially an equidistant mesh with mesh size $h$ in $[0,\tau]$ and
with mesh size $H$ in $[\tau,1]$.  We choose for $k$-th order elements $\tau=\eps^{(k-1)/k}$.The fine mesh is given by
\[
x_i=ih,\,\,i=0,1,\cdots,[1/H]\quad{\rm with}\quad  h=\alpha H\eps^{(k-1)/k}
\]
This choice works
works because we obtain on the fine mesh
\[
|E-E^I|_1\pr h^k|E|_{k+1}\pr \eps^{k-1}H^k\eps^{-(k-1/2)}\pr \eps^{-1/2}H^k
\]
and
\[
|E-E^I|_0\pr h^k|E|_{k}\pr \eps^{1/2}H^k.
\]
If the layer part $E$ in the transition point $\tau$ is sufficiently small we get
\be
 \|u-u_h\|_\eps \pr H^k.
\ee
Let us assume $k=2$. Then the smallness condition is satisfied if
\[
\eps\pr H^3 \quad {\rm or}\quad \eps e^{-1/(\eps^{1/2})}\pr H^3.
\]
\begin{remark}
The second condition is not very restrictive. For instance, the condition reads
\[
5\times 10^{-12}\le H^3 \quad {\rm if}\,\,\eps\le 0.0025.
\]
If $k$ increases, this condition becomes more and more
restrictive. This means, our approach makes sense if $\eps$ is extremely small or $k$ is only of
moderate size.
\end{remark}

For the reaction-diffusion problem, we get analogously
\be
 \|u-u_h\|_\eps \pr \eps^{1/2}H^k+H^{k+1}.
\ee
Using the approach of \cite{FR20} it is also possible to prove some result in the balanced norm.
\begin{remark}
Our mesh is not {\it locally uniform}. To get this property, instead of the uniform mesh in $[\tau,1]$
one could use a graded mesh with
\[
  x_{i+1}=(1+H)x_i\quad{\rm for}\quad i\ge [1/H],
\]
following \cite{DL06}. But then the number of mesh points used depends on $\ln(1/\eps)$.
\end{remark}

\section{Fourth order problems and cubic $C^1$-splines}
In many fourth order problems typically weak, but no strong layers exist. For finite element methods
on layer adapted meshes, see \cite{FR20,SSa95,SSb95,Xe17}.\\
We start with a problem of {\it convection-diffusion} type:
\be
\eps u^{(4)}+bu'''+L_2u=f, \quad u(0)=u''(0)=u(1)=u''(1)=0.
\ee
Here $L_2$ is a linear second order operator and we assume that the bilinear form associated to the full
operator is V-elliptic and $b>1$. Then, it is well known that we have a solution decomposition into a
smooth part and a layer part with (see \cite{Ga88,SSa95,SSb95}).
\[
    |E^{(k)}|\preceq \eps^{2-k}e^{-x/\eps}.
\]
That means related to the given boundary conditions the layer is very weak. Now we use the norm
\[
   \|v\|\eps:=\eps^{1/2}|v|_2+|v|_1.
\]
If we choose $\tau=\eps^{1/2}$ and $h=\alpha\eps^{1/2}H$, we get for the interpolation error $\eta$ on
the fine mesh the estimates
\[
  |\eta|_2\pr h^2|E|_4\pr h^2\eps^{-3/2}\pr \eps^{-1/2}H^2
\]
and
\[
  |\eta|_1\pr h^2|E|_3\pr h^2\eps^{-1/2}\pr \eps^{1/2}H^2.
\]
These estimates allow us to prove
\be
 \|u-u_h\|_\eps \pr H^2.
\ee
\begin{remark}
If the given equation is equipped with such boundary conditions that the layer is extremely weak with
\[
    |E^{(k)}|\preceq \eps^{3-k}e^{-x/\eps},
\]
we have for the interpolation error on an {\it equidistant} mesh
\[
|\eta|_2\pr \eps^{-1/2}H^2\quad {\rm and}\,\, |\eta|_1\pr \eps^{1/2}H^2.
\]
This allows to prove uniform convergence if the complete boundary conditions allow
to prove V-ellipticity of the related bilinear form.
\end{remark}
But if we next study the problem
\be
\eps u^{(4)}+bu'''+L_2u=f, \quad u(0)=u'(0)=u(1)=u'(1)=0.
\ee
with a weak layer, i.e., with
\[
    |E^{(k)}|\preceq \eps^{1-k}e^{-x/\eps},
\]
the situation becomes different. We have for the interpolation error of $E$ on the fine mesh
\[
   \eps^{1/2}|\eta|_2\pr \eps^{1/2}h^2|E|_4\pr h^2\eps^{-2}.
\]
That means, to achieve a second order result, one should use $h=\alpha\eps H$. This leads to $\tau=\eps$.
But then $E$ is in $\tau$ only small enough if $\eps$ is extremely small. Consequently, for that problem
one should prefer a Shishkin type mesh with $\tau=\tau_0\eps\ln(1/H)$.

Consider finally problems of {\it reaction-diffusion} type
\be
\eps^2 u^{(4)}+L_2u=f, \quad u(0)=u''(0)=u(1)=u''(1)=0.
\ee
Assume ellipticity of $L_2$. Now we have two layers at $x=0$ and at $x=1$, and assume, for instance for the
layer E at $x=0$
\[
    |E^{(k)}|\preceq \eps^{2-k}e^{-x/\eps}.
\]
Because now we estimate in the norm
\[
   \|v\|\eps:=\eps|v|_2+|v|_1,
\]
we get with $\tau=\eps^{1/2}$
\be
 \|u-u_h\|_\eps \pr \eps^{1/2}H^2+H^3.
\ee
Using \cite{FR20}, an estimate in the balanced norm should also be possible.

In the case of a weak layer of the boundary value problem
\be
\eps^2 u^{(4)}+L_2u=f, \quad u(0)=u'(0)=u(1)=u'(1)=0
\ee
we have
\[
   \eps|\eta|_2\pr \eps h^2|E|_4\pr h^2\eps^{-3/2}.
\]
Consequently,
we get for the choice $\tau=\eps^{3/4}$ the estimate
\be
 \|u-u_h\|_\eps \pr H^2.
\ee
But we see no possibility to prove a balanced norm estimate.

\section{A mixed finite element method for some fourth order problems}

Similarly as in \cite{FR16}, we consider a mixed method for the problem
\be
\eps^2 u^{(4)}-bu''+du=f, \quad u(0)=u'(0)=u(1)=u'(1)=0.
\ee
We assume
\[
d-\frac{1}{2}b''>\delta>0,
\]
moreover, the existence of a decomposition
\be
u=S+E_1+E_2 \quad {\rm with}\quad |E_1^{(k)}|\pr \eps^{1-k}e^{-x/\eps}
\ee
and the corresponding estimate for the layer $E_2$ at $x=1$.

Introducing $w=\eps u''$, the mixed method is based on:\\
Find $(u,w)\in H_0^1\times H^1$ such that
\be
 \eps(u',\phi')+(w,\phi)=0 \quad \forall \phi\in H^1
\ee
and
\be
(bu',\psi')+(du,\psi)-\eps(w',\psi')=(f,\psi) \quad \forall \psi\in H_0^1.
\ee
We have coercivity in the norm
\be\label{norm}
 \|(u,w)\|_1^2:=|u|_1^2+\|w\|_0^2.
\ee
$u$ has weak layers, the layers of $w$ are strong. But because in the norm \eqref{norm} only
the $L_2$ norm of $w$ appears, we expect that it is possible to use a mesh coarser than a standard
Shishkin type mesh used in \cite{FR16}.

Consequently, we consider a mixed finite element method with $P_k$-elements for $u$ and $w$ to obtain the
discrete solution
$(u_H,w_H)$ on a special mesh. The mesh is equidistant and fine in $[0,\tau]$ and $[1-\tau,1]$, in the
remaining part the mesh is equidistant with the mesh size $H$.

To define $\tau$, we study first the {\it interpolation error}.
 It is sufficient to consider the error generated by $E_1$ on $[0,\tau]$, because on the remaining part
 of the interval $E_1$ is sufficiently small and therefore the $L_\infty$ stability of the interpolation
 operator guarantees the same for the interpolation error.

Similarly as in Section 2 we get for $\tau=\eps^{1-1/(2k)}$ and $h=\alpha H\eps^{1-1/(2k)}$ for the interpolation
error
\be\label{int}
\|(u-u^I,w-w^I)\|_1 \pr H^k.
\ee
The smallness of the layers in the transition points makes our approach useful for $k=1,2$.

For the study of the {\it discrete error} $(\psi_H,\phi_H)=(\pi u-u_H,\pi w-w_H)$ we first let the choice of
the interpolation operator $\pi$ open. Of course, we assume that $\pi$ has the same approximation error
properties as the Lagrange interpolation and, consequently, \eqref{int} holds as well for $u-\pi u,w-\pi w$.
For the discrete error we obtain with $\eta=\pi u-u$ and $\zeta=\pi w-w$ (see \cite{FR16})
\be\label{error}
\|(\psi_H,\phi_H)\|_1^2\pr \eps(\eta',\phi'_H)+(\zeta,\phi_H)+(b\eta',\psi'_H)
   +(d\eta,\psi_H)-\eps(\zeta',\psi'_H).
\ee
For linear elements we have
\[
  ((u-u^I)',\phi'_H)=0,
\]
because $\phi'_H$ is piecewise constant and $u-u^I$ vanishes in all mesh points.
This property simplifies the error estimation.\\
 Therefore, we introduce for $k\ge 2$ a new interpolant $\pi$.
This interpolant satisfies on the interval $[x_{i-1},x_i]$ first
$(\pi v)(x_{i-1})=v(x_{i-1})$ and $(\pi v)(x_{i})=v(x_{i})$, moreover
\[
  \int_{x_{i-1}}^{x_i}(x-x_{i-1})^l \pi v=\int_{x_{i-1}}^{x_i}(x-x_{i-1})^l  v \quad {\rm for}\,\,
      l=1.\cdots,k-1.
\]
Then, the first and the last term in \eqref{error} vanish, see Lemma 2.66 in \cite{RST08}..
 Moreover, the interpolant
has the standard approximation properties and is $L_\infty$ stable. Equation \eqref{error} reduces to
\be
\|(\psi_H,\phi_H)\|_1^2\pr (\zeta,\phi_H)+(b\eta',\psi'_H)
   +(d\eta,\psi_H),
\ee
and with
\[
   |\eta|_1\pr H^k \quad {\rm and}\quad \|\zeta\|_0\pr H^k
\]
one gets easily
\be\label{fin}
\|(\pi u-u_H,\pi w-w_H)\|_1 \pr H^k\quad {\rm and}\quad \|(u-u_H,w-w_H)\|_1 \pr H^k.
\ee

Analogously one can handle the case of very weak layers of the problem
\be
\eps^2 u^{(4)}-bu''+du=f, \quad u(0)=u''(0)=u(1)=u''(1)=0.
\ee
Then, for linear elements we derive on a uniform mesh
\be
\|(u-u_H,w-w_H)\|_1 \pr H.
\ee
If $k\ge 2$, we choose $\tau=\eps^{1-3/(2k)}$ and obtain again \eqref{fin}.

\end{document}